\documentclass[a4paper,reqno]{amsart}
\usepackage{amssymb}

\newtheorem{theorem}{Theorem}
\newtheorem{lemma}{Lemma}
\newtheorem{cor}{Corollary}
\newtheorem{prop}{Proposition}
\theoremstyle{remark}
\newtheorem*{ack}{Acknowledgements}
\newtheorem*{rem*}{Remark}
\newtheorem*{ex*}{Example}

\newtheoremstyle{citing}
  {3pt}
  {3pt}
  {\itshape}
  {}
  {\bfseries}
  {.}
  {.5em}
  {\thmnote{#3}}
\theoremstyle{citing}
\newtheorem*{varthm}{}

\newcommand{\set}[1]{\{#1\}}
\newcommand{\sett}[2]{\{#1\,|\,#2\}}
\DeclareMathOperator{\l<}{\langle}
\DeclareMathOperator{\r>}{\rangle}
\newcommand{\gen}[1]{\l<#1\r>}
\newcommand{\field}[1]{\mathbb{#1}}
\newcommand{\Q}{\field{Q}}
\newcommand{\C}{\field{C}}
\newcommand{\N}{\field{N}}
\newcommand{\Z}{\field{Z}}
\newcommand{\F}{\field{F}}
\newcommand{\bs}{\backslash}
\newcommand{\ift}{\infty}
\newcommand{\ndiv}{\nmid}
\newcommand{\lra}{\longrightarrow}
\newcommand{\lms}{\longmapsto}
\newcommand{\ovl}[1]{\overline{#1}}
\newcommand{\wt}{\widetilde}

\newcommand{\bt}{{\mathbf t}}

\DeclareMathOperator{\Hom}{Hom}
\DeclareMathOperator{\Aut}{Aut}
\DeclareMathOperator{\Gal}{Gal}
\DeclareMathOperator{\charak}{char}
\DeclareMathOperator{\ord}{ord}
\DeclareMathOperator{\im}{im}

\newcommand{\gwgwp}{generalized weak Grunwald-Wang problem }
\newcommand{\gwgwpp}{generalized weak Grunwald-Wang problem}

\title{On absolute Galois splitting fields of central simple algebras}

\author{Timo Hanke${}^1$}
\address{
Department of Mathematics\\
Technion -- Israel Institute of Technology\\
Haifa, 32000\\
Israel.
email : hanke@math.uni-potsdam.de
}
\date{\today}

\begin{document}

\begin{abstract}
  A splitting field of a central simple algebra is said to be {\em absolute Galois}
if it is Galois over some fixed subfield of the centre of the algebra.
The paper provides an existence theorem for such fields over global fields
with enough roots of unity.
As an application, all twisted function fields and all twisted Laurent series rings over symbol algebras (or $p$-algebras) 
over global fields are crossed products.
A closely related statement holds for
division algebras over Henselian valued fields with global residue field.

The existence of absolute Galois splitting fields in central simple algebras over global fields
is equivalent to a suitable generalization of the weak Grunwald-Wang Theorem,
which is proved to hold if enough roots of unity are present.
In general, it does not hold and counter examples have been used in noncrossed product constructions.
This paper shows in particular 
that a certain computational difficulty involved in the construction of explicit examples of noncrossed product twisted Laurent series rings
can not be avoided by starting the construction with a symbol algebra.
\end{abstract}

\maketitle
\begin{quotation}
  \tiny NOTICE: this is the author's version of a work that was accepted for
publication in Journal of Number Theory. Changes resulting from the publishing
process, such as peer review, editing, corrections, structural formatting, and
other quality control mechanisms may not be reflected in this document. Changes
may have been made to this work since it was submitted for publication. 
A definitive version will be subsequently published in Journal of Number Theory (2007), doi:10.1016/j.jnt.2006.12.011.
\end{quotation}

\footnotetext[1]{Supported by the Technion -- Israel Institute of Technology.}
\section*{Introduction}

A global field is a finite extension of $\Q$ or $\F_p(\bt)$
and by a valuation on it we mean a nontrivial
absolute value in the nonnegative real numbers,
archimedian or non-archimedian.
Let $L/K$ be an extension of global fields. 
If $v$ and $w$ are valuations on $K$ and $L$, respectively,
then we write $w|v$ if $w$ extends $v$. 
In this case $L_w/K_v$ denotes the corresponding extension of completions
and $[L:K]_w$ denotes the local degree $[L_w:K_v]$. 
Clearly,
\begin{equation}
 [L:K]_w \leq\begin{cases}
  [L:K] \\ (2,[L:K]) \\ 1
\end{cases}
\text{if $v$ is }
\begin{cases}
  \text{non-archimedian,}\\
  \text{real archimedian,}\\
  \text{complex archimedian,}
\end{cases}
  \label{eq:bound}
\end{equation}
where parentheses stand for the greatest common divisor.
We say that $L/K$ has {\em full local degree at $v$} 
if the corresponding bound in \eqref{eq:bound} is achieved for all $w|v$.
For non-archimedian $v$ this is equivalent to saying that $v$ extends uniquely to $L$.
Throughout the paper $S$ denotes a finite set of pairwise inequivalent valuations on $K$.
The extension $L/K$ is said to have {\em full local degree in $S$} 
if $L/K$ has full local degree at each $v\in S$.
We refer to the following well-known statement as the

\begin{varthm}[Weak Grunwald-Wang Theorem]
  \label{thm:wGW}
For any finite $S$ and any $m\in\N$
there is a cyclic extension $L/K$ of degree $m$ with full local degree in $S$.
\end{varthm}

\begin{rem*}
The original theorem of Grunwald-Wang prescribes not only local degrees
but also the local completions $L_v$ up to isomorphism at each $v\in S$.
Most references on the weak version
(e.g.\ 
Hasse \cite[Ganz schwacher Existenzsatz]{hasse:grunwald} or 
Artin-Tate \cite[Ch.\ 10, Thm.\ ~5]{artin-tate})
state a form of the theorem
that can prescribe any locally possible degrees instead of just full degrees.
Note that the weak version has no special case.
\end{rem*}

Recall that the {\em degree} of a finite-dimensional central simple $K$-algebra $A$,
denoted $\deg A$,
is the square root of the dimension of $A$ over $K$.
A {\em splitting field} of $A$ 
is a field $L\supseteq K$ such that $A\otimes_K L\cong M_m(L)$ where $m=\deg A$.
A {\em maximal subfield}\footnote{maximal as a commutative subring of $A$} of $A$ is a field $L\subseteq A$
with $[L:K]=\deg A$. 
Any maximal subfield is a splitting field
and any splitting field with $[L:K]=\deg A$ is a maximal subfield.
According to the local-global principle for central simple algebras over global fields 
(c.f.\ Pierce \cite[\S18]{pierce:ass-alg})~:
\begin{equation*}
    \begin{split}
  &\text{$L/K$ has full local degree in $S \iff$}\\
  &\text{$L$ is a splitting field of all central simple $K$-algebras $A$ with $\deg A=[L:K]$}\\
  &\text{and Hasse invariant zero outside of $S$.}
    \end{split}
\end{equation*}
The set of nonzero Hasse invariants is always finite.
Consequently, 
every central simple $K$-algebra is cyclic
(i.e.\ has a maximal subfield cyclic over $K$).
This statement is considered the primary application of the weak Grunwald-Wang Theorem.

Now, let $k\subseteq K$ be a subfield such that $K/k$ is Galois
and let $m\in\N$.
As before, $S$ is a finite set of pairwise inequivalent valuations on $K$.
We denote by $(K/k,m,S)$ the 
\begin{varthm}[Generalized Weak Grunwald-Wang Problem]
  \label{problem}
  Does there exist a cyclic extension $L/K$ of degree $m$ 
  with full local degree in $S$ such that $L/k$ is Galois ?
  Any such $L$ is said to be a solution to $(K/k,m,S)$.
  A solution is called {\em abelian (cyclic)} if $L/k$ is abelian (cyclic).
\end{varthm}

\begin{ex*}
$(\Q(\zeta_7+\zeta_7^{-1})/\Q,3,\set{2})$ has no solution.\footnote{
Here, $\zeta_7$ denotes a primitive $7$-th root of unity.
We briefly sketch the argument~:
if $L$ is a solution then $L/k$ is either cyclic or elementary abelian
because the degree $[L:K]$ is $9$.
There is no cyclic solution because the prime $7$ ramifies in $\Q(\zeta_7+\zeta_7^{-1})$ 
and $\Q_7$ has no totally ramified cyclic extension of degree $9$.
The prime $2$, on the other hand, is inert in $\Q(\zeta_7+\zeta_7^{-1})$.
Since $\Q_2$ has no totally ramified cyclic cubic extension,
$L/k$ can not be elementary abelian and have full local degree at the $2$-adic valuation.
}
\end{ex*}

We say that a splitting field or maximal subfield of a central simple $K$-algebra $A$ is
{\em $k$-Galois} if it is Galois over $k$.
The term {\em absolute Galois} is used alternatively 
if there is no ambiguity about $k$.
Recall that $A$ is called a {\em crossed product} if $A$ contains a $K$-Galois maximal subfield.
The existence of a $k$-Galois maximal subfield is thus stronger than $A$ being a crossed product.
According to the local-global principle, 
$(K/k,m,S)$ has a solution for all $S$ if and only if every central simple $K$-algebra of degree $m$
has a $k$-Galois splitting field of degree $m$ over $K$,
if and only if every central simple $K$-algebra contains a $k$-Galois maximal subfield.

Absolute Galois splitting fields occur naturally
in the classification of crossed products over valued fields.
For instance, if $k$ is an arbitrary perfect field 
then any division algebra central over $k((\bt))$ 
is a crossed product if and only if 
the residue division algebra (with respect to the $\bt$-adic valuation) 
contains a $k$-Galois maximal subfield.
As another example, 
any twisted Laurent series ring $A((\bf x;\sigma))$
or twisted function field $A(\bf x;\sigma)$ over a simple algebra $A$
with $\sigma\in\Aut(A)$
is a crossed product if and only if $A$ contains a $k$-Galois maximal subfield
where $k$ is the field of central invariants $Z(A)^\sigma$.
These two examples are corollaries from the following theorem due to Hanke \cite[Thm.\ 5.20]{hanke:thesis}. 
\begin{theorem}
  Let $F$ be a Henselian valued\footnote{The valuation in this theorem can be an arbitrary Krull valuation.} field
  with residue field $\ovl F$.
  Any $F$-central division algebra $D$
  that splits over the maximal unramified extension $F_{nr}$
  is a crossed product if and only if the residue division algebra $\ovl D$
  contains an $\ovl F$-Galois maximal subfield.
  \label{thm:0}
\end{theorem}
In the situation of Theorem \ref{thm:0}, 
the centre $Z(\ovl D)$ is known to be an abelian extension of $\ovl F$
with the abelian rank bounded by the rank of the valuation
(cf.\ \cite[Lem.\ 5.1]{jacob-wadsworth:div-alg-hensel-field}).
In particular, $Z(\ovl D)/\ovl F$ is cyclic if the valuation is discrete
(for instance, in the two examples given before Theorem \ref{thm:0}).
This creates a special interest in the case $K/k$ cyclic of the \gwgwpp.
It shall be pointed out that despite the connection made by Theorem \ref{thm:0},
the present paper
does not require any knowledge about valuations on division rings
(with the exception of the last corollary).

Consequently from Theorem \ref{thm:0}, division algebras without absolute Galois maximal subfield
can lead to noncrossed products.
More precisely, the case $K/k$ cyclic can give noncrossed products over $k(\bt)$ and $k((\bt))$,
while the case $K/k$ abelian of rank $r$ 
can give noncrossed products over $k(\bt_1,\ldots,\bt_r)$ and $k((\bt_1,\ldots,\bt_r))$.
This was first demonstrated by E.\ Brussel \cite{brussel:noncr-prod}
who proved for number fields $k$ the existence of cyclic extensions $K/k$
and $K$-central division algebras $D$ without $k$-Galois maximal subfield
(cf.\ \cite[Lem.\ 5]{brussel:noncr-prod})
and derived from it the existence of noncrossed products over $k(\bt)$ and $k( (\bt) )$.
The realizable range of degrees of $K/k$ and $D/K$ depends on the number of roots of unity present in $k$.
Another example with $K/k$ bicyclic and $D$ a quaternion algebra can be found in Hanke \cite{hanke:expl-ex}.

A certain type of $k$-Galois splitting fields,
namely those that are composites of $K$ with a cyclic extension of $k$,
were studied already by Albert \cite[\S 9]{albert:p-adic-fields-and-rational-division-algebras}.
Later, Liedahl \cite{liedahl:qi} characterized the $\Q(i)$-division algebras with $\Q$-Galois maximal subfields.
Note that an error in Albert's assumptions was discovered by Liedahl
and the statement was corrected in \cite[Prop.\ 6]{liedahl:qi}.

This paper is concerned exclusively with the case $K/k$ cyclic of the \gwgwp and proves the following theorem. 
The notation $\mu_m\subset F$ means that a field $F$ contains all $m$-th roots of unity of its algebraic closure.

\begin{varthm}[Main Theorem]
 Let $K/k$ be a cyclic extension of global fields and let $m\in\N$. 
 For any $S$,
 \begin{enumerate}
   \item 
 if $\mu_m\subset K$ then $(K/k,m,S)$ has a solution.
 \end{enumerate}
 Suppose $\mu_m\subset k$. For any $S$,
 \begin{enumerate}
   \setcounter{enumi}{1}
   \item $(K/k,m,S)$ has an abelian solution,
   \item $(K/k,m,S)$ has a cyclic solution if and only if $\mu_m\subset N_{K/k}(K)$.
 \end{enumerate}
 In particular, $(K/k,m,S)$ has a cyclic solution if $\charak k=p$ and $m$ is a $p$-power. 
  \label{thm:main}
\end{varthm}

The Main Theorem can be regarded as the explanation why Brussel's proof 
for the existence of noncrossed products in \cite{brussel:noncr-prod} 
depends on the absence of certain roots of unity.
The assumption on $K/k$ being cyclic is necessary.
Indeed,
it is shown in \cite[Thm.\ 6.2]{hanke:expl-ex} that 
$(\Q(\sqrt{3},\sqrt{-7})/\Q,2,\set{3,7})$ has no solution.

Because of the presence of roots of unity
the proof of the Main Theorem requires only Kummer theory (resp.\ Witt vectors) 
plus the approximation theorem.
This is in analogy to the Grunwald-Wang Theorem
which also has such an elementary proof if enough roots of unity 
are contained in the ground field
(cf.\ Hasse \cite[\S 4]{hasse:grunwald}).
As a consequence the Main theorem holds in principle over arbitrary multi-valued fields instead of global fields
(in the sense of Lorenz-Roquette \cite{lorenz-roquette:gw},
cf.\ the remark made there at the end of \S 1.1).
Only, one has to assume the existence of enough cyclic extensions over the completions
as it is made precise in Theorem \ref{thm:exist-a}.
However, this kind of generalization will not be pursued in this paper.

For the following application recall that a {\em symbol algebra} of degree $m$
over a field $K$ with $(\charak K,m)=1$ and $\mu_m\subset K$
is a central simple $K$-algebra generated by two elements $x,y$ subject to the relations
$$ x^m=a,\quad y^m=b,\quad yx=\zeta xy$$
for some $a,b\in K^*$ and $\zeta$ a primitive $m$-th root of $1$.
A {\em $p$-algebra}, $p$ prime, is a central simple algebra of $p$-power degree 
over a field of characteristic $p$.

\begin{varthm}[Corollary]
  Let $A$ be a central simple algebra over a global field.
  If $A$ is a tensor product of a symbol algebra with a $p$-algebra
  (the factors are allowed to be trivial) then 
  any twisted Laurent series ring $A( (\bf x;\sigma) )$ and 
  any twisted function field $A(\bf x;\sigma)$  over $A$ is a crossed product.
\end{varthm}

The corollary is a complementary result 
to the noncrossed product twisted Laurent series ring $D((\bf x;\sigma))$ from Hanke \cite{hanke:laurent-noncr},
where $D$ is a cyclic cubic division algebra.

\begin{rem*}
Wang \cite{wang:gw} discusses cyclic solutions to the
``strong'' generalized Grunwald-Wang problem,
where not only local degrees but also the local completions are prescribed.
\end{rem*}

The structure of the paper is as follows. 
Section $1$ provides an existence statement (Theorem \ref{thm:exist-a})
for certain local elements that will serve as Kummer radicands in our construction; 
section $2$ proves the Main Theorem split up into 
the statement on cyclic solutions (Theorem \ref{thm:cyclic}) 
and the statement on general and abelian solutions (Theorem \ref{thm:main2});
section $3$ formulates the Main Theorem in terms of $k$-Galois splitting fields of central simple algebras (Theorem \ref{thm:main-cs})
and derives two versions of a crossed product characterization 
(Corollaries~\ref{cor:symbol2} and \ref{cor:hensel}).
There is an appendix on Albert's Theorem on the cyclic embedding problem.

\begin{ack}
I would like to thank the Department of Mathematics of the Technion -- Israel Institute of Technology, Haifa,
for their kind hospitality and financial support during my visit in 2005.
Particular thanks go to Jack Sonn for the many helpful discussions.
Furthermore, I am indebted to the referee for suggesting
improvements of the exposition.
\end{ack}

\section{Preliminaries}

\subsection{Notations}

If $k$ is a field then $k^*$ is the multiplicative group of nonzero elements.
For any $m\in\N$ we denote by $\mu_m(k)$ the group of all $m$-th roots of unity contained in $k$.
The statement $\mu_m\subset k$ is short for $\mu_m(\ovl k)\subset k$ where $\ovl k$ is the algebraic closure of $k$.
Note that we include in this notation the case $\charak k|m$
in which $\mu_m\subset k$ is equivalent to $\mu_{m_0}\subset k$
where $m_0$ is the greatest divisor of $m$ prime to $\charak k$.

Unless stated otherwise,
a valuation is meant to be a nontrivial archimedian or non-archimedian
valuation of rank 1,
i.e. the values are real numbers.
Throughout the paper $S$ denotes a finite set of pairwise inequivalent
valuations on the global field $K$.
The complex archimedian valuations $v$ on $K$
are not of interest for the matter of splitting fields
because all central simple $K$-algebras are already split over the completion $K_v=\C$.
We therefore assume for the sake of simplicity 
that no such valuations occur in $S$
(in fact, a few statements of this paper would require slight modifications if complex valuations were allowed in $S$).
Define $S_0\cup S_\ift$ to be the partition of $S$
into its non-archimedian and archimedian valuations.
If $v$ is a valuation on $K\supseteq k$ we simply denote the restriction of $v$ to $k$ also by $v$.
The completions with respect to $v$ are denoted $k_v$ and $K_v$, respectively.

\subsection{Irreducible Binomials}

Let $k$ be field and consider the binomials $x^n-a$ with $a\in k$.
Let $p$ denote a prime unless stated otherwise.
We recall a well-known theorem 
(e.g.\ Lang~\cite[Ch.\ VIII, \S 9]{lang:algebra})~: 
\begin{equation}
  x^n-a \text{ irreducible over $k \iff a\not\in k^p$ for all $p|n$ and $a\not\in -4k^4$ if $4|n$.}
  \label{eq:standard}
\end{equation}
Obviously, $a\not\in -k^2$ implies $a\not\in -4k^4$.
The condition $a\not\in -k^2$ is not only more convenient to check,
it also carries over to Galois extensions of $k$ in the following sense.

\begin{prop}
  \label{prop:forall-gal}
  Let $a\in k$ and $n\in\N$.
  If $4$ divides $n$ we assume $a\not\in-k^2$.
  Then one has
  \begin{equation}
    \begin{split}
    &\text{for all Galois extensions $K/k$ :}\\
    &x^n-a \text{ irreducible over $K$} \iff a\not\in K^p \text{ for all $p|n$.}
    \end{split}
    \label{eq:star}
  \end{equation}
\end{prop}
Before starting with the proof we make a few general observations. 
As an immediate consequence of \eqref{eq:standard},
  \begin{equation}
    \begin{split}
&x^n-a \text{ irreducible over $k$} \iff\\
&x^p-a \text{ irreducible over $k$ for all $p|n$ with $p$ prime or $p=4$.}
    \end{split}
  \label{eq:p-4}
  \end{equation}
If $\charak k\neq 2$ then trivially
\begin{equation}
  a\not\in-4k^4 \iff x^4+4a \text{ has no root in $k$,}
  \label{eq:x4-no-root}
\end{equation}
and by \eqref{eq:standard},
\begin{equation}
  a\not\in-k^2 \text{ and } a\not\in k^4 \iff x^4+4a \text{ irreducible over $k$.}
  \label{eq:x4-irred}
\end{equation}

\begin{proof}[Proof of Proposition \ref{prop:forall-gal}]
  The case $\charak k=2$ or $4\ndiv n$ of the assertion is already included in \eqref{eq:standard},
  so we assume throughout the proof $\charak k\neq 2$ and $4$ divides $n$.
  Hence $a\not\in -k^2$ is also assumed by hypothesis.
  Let $K/k$ be a Galois extension.
  Owing to \eqref{eq:standard}, it remains to show that $a\not\in K^2$ implies $a\not\in-4K^4$.
  Let $a\not\in K^2$. 
  In particular, $a\not\in k^4$ and $x^4+4a$ is irreducible over $k$ by \eqref{eq:x4-irred}. 
  If $a\in-4K^4$ we get a contradiction as follows.
  On the one hand, $a\in-K^2$, hence $-1\not\in K^2$.  
  On the other hand, since the irreducible polynomial $x^4+4a$ has a root in $K$
  and $K/k$ is Galois, it splits completely over $K$.
  Hence $-1\in K^2$, a contradiction.
\end{proof}

The following remark is not further used in this paper but stated for completeness.
It shows that the assumption of Proposition \ref{prop:forall-gal} is necessary.

\begin{rem*}
  Let $a\in k$ and $n\in\N$.
  If $a\not\in k^2$ and $4$ divides $n$ then \eqref{eq:star} implies $a\not\in -k^2$.
\end{rem*}
\begin{proof}
  Suppose $a\not\in k^2, 4|n$ and \eqref{eq:star} holds.
  Assume $\charak k\neq 2$,
  for otherwise $a\not\in-k^2$ is trivial.
  We assume $a=-b^2$ for some $b\in k$ and derive a contradiction.
  The fields $K_{1,2}=k(\sqrt{\pm 2b})$ both contain a root of $x^4+4a$.
  Hence, $a\in -4K_i^4$ by \eqref{eq:x4-no-root}.
  Since $4|n$ this means by \eqref{eq:standard} that $x^n-a$ is reducible over $K_i$.
  The property \eqref{eq:star} then implies $a\in K_i^2$.
  Since $a\not\in k^2$ we have $K_i=k(\sqrt{a})$,
  in particular $K_1=K_2$.
  By construction of $K_{1,2}$ it follows $-1\in k^2$.
  Since $a\in-k^2$ we get $a\in k^2$
  in contradiction to the hypothesis.
\end{proof}

\subsection{Existence of radicands}

\begin{theorem}
  Let $k$ be a field and let $p$ be a prime different from $\charak k$.
  We assume $k\neq k^p$ and if $\mu_p\subset k$ we assume $k^*/{k^*}^p$ not cyclic.
  For any cyclic extension $K/k$ and any $s\in\N$
  there is an $a\in k^*$ such that the binomial $x^{p^s}-a$ is irreducible over~$K$.
  \label{thm:exist-a}
\end{theorem}
\begin{proof}
  First observe that the hypothesis on $k$ implies $U:={k^*}^p\gen{-1}\neq k^*$.
  For, if $p\neq 2$ then $U={k^*}^p\neq k^*$,
  and if $p=2$ then $k^*/{k^*}^p$ is not cyclic.
  
  Let $K/k$ be any cyclic extension and define $V:=K^p\cap k^*$.
  We show that $V\neq k^*$.
  If $V=k^*\neq {k^*}^p$ then 
  the Galois extension $K/k$ contains a Kummer subfield $K_0$ of degree $p$ over $k$.
  In particular, $\mu_p\subset k$.
  Moreover, since $K_0$ is the unique subfield of degree $p$, we have $K_0=k(V^{1/p})$.
  Hence, $V/{k^*}^p$ is cyclic by Kummer theory, contrary to the hypothesis.

  Since $U$ and $V$ are proper subgroups of $k^*$
  we have $U\cup V\neq k^*$.
  This proves the Proposition.
  For, if $a\in k^*\bs(U\cup V)$ then $x^{p^s}-a$ is irreducible over $K$ by Proposition~\ref{prop:forall-gal}.
\end{proof}

\begin{rem*}
  The assumptions of Theorem \ref{thm:exist-a} are necessary.
  Indeed, if $\mu_p\subset k$ and $k^*/{k^*}^p$ is cyclic
  then by Kummer theory 
  there is a unique cyclic extension $K/k$ of degree $p$.
  Hence $x^p-a$ with $a\in k^*$ is never irreducible over that $K$.

  Also, Theorem~\ref{thm:exist-a} does not generalize to non-cyclic extensions $K/k$.
  For instance, the field $\Q_3$ of $3$-adic numbers has exactly three quadratic extensions.
  Thus, Theorem~\ref{thm:exist-a} does not hold for $k=\Q_3$, $p=2$ and $K$ the
unique biquadratic extension of $k$.
\end{rem*}

\begin{lemma}
  Let $k$ be a field with a non-archimedian discrete valuation.
  For any prime $p$ we have $k\neq k^p$ and 
  \begin{equation}
    \text{$k^*/{k^*}^p$ is cyclic} \iff U_k\subseteq {k^*}^p,
    \label{eq:discr-val-cond}
  \end{equation}
  where $U_k$ denotes the valuation units.
  \label{lem:discr-val-cond}
\end{lemma}
\begin{proof}
  Let $a\in k^*$ be a uniformizer.
  Plainly, $k\neq k^p$ because $a\not\in k^p$.
  Suppose $U_k\not\subseteq {k^*}^p$ and let $u\in U_k, u\not\in{k^{*}}^p$.
  Then $k^*/{k^*}^p$ is not cyclic because it is a group of exponent $p$ and $a\not\in {k^*}^p\gen{u}$.
  Conversely, suppose $U_k\subseteq {k^*}^p$.
  Since for any $x\in k^*$ we have $x=a^ru$ with $r\in\Z$ and $u\in  U_k$
  it follows $k^*={k^*}^p\gen{a}$.
\end{proof}

A local field is the completion of a global field with respect to some non-archimedian valuation,
i.e.\ a finite extension of $\Q_l$ or $\F_l((\bt))$.

\begin{cor}
  \label{cor:ass}
  Any local field satisfies the assumptions of Theorem \ref{thm:exist-a}.
\end{cor}
\begin{proof}
  Let $k$ be a local field and $p\neq\charak k$.
  If $U_k\subseteq {k^*}^p$ and $\mu_p\subset k$ then $k$ contains all $p$-power roots of $1$.
  Since this is obviously not the case for finite extensions of $\Q_l$ or $\F_l((\bt))$,
  the statement follows from Lemma \ref{lem:discr-val-cond} (the valuation is discrete). 
\end{proof}

\section{Proof of the Main Theorem}

We divide the Main Theorem into the statement about cyclic solutions
(Theorem \ref{thm:cyclic})
and the statements about general and abelian solutions
(Theorem \ref{thm:main2}),
reflecting the different nature of their proofs.
Moreover, the proof for general and abelian solutions makes use of the cyclic result
when dealing with the degree part that is a power of $\charak k$.

\subsection{Cyclic solutions}

Let $K/k$ be a cyclic extension of global fields and let $m\in\N$. 

\begin{theorem}
  For any $S$,
  if $\mu_m\subset k$ then $(K/k,m,S)$ has a cyclic solution if and only if $\mu_m\subset N_{K/k}(K)$.
  In particular, there is a cyclic solution if $m$ is a power of $\charak k$.
  \label{thm:cyclic}
\end{theorem}

The rest of this section is devoted to the proof of Theorem \ref{thm:cyclic}.
We start by giving two easy lemmas which reduce consideration 
to the case where $[K:k]$ and $m$ are both powers of the same prime number.
For any prime $p$ we denote by $m_p$ the maximal $p$-power dividing $m$,
by $K_p$ the unique subfield of $K/k$ with maximal $p$-power degree,
and by $S_p$ the finite set of pairwise inequivalent valuations on $K_p$
which is obtained from $S$ by restricting each valuation to $K_p$.

\begin{lemma}
  $(K/k,m,S)$ has a cyclic solution if and only if $(K_p/k,m_p,S_p)$ has a cyclic solution for each prime divisor $p$ of $m$.
  \label{lem:p}
\end{lemma}
\begin{proof}
  If $L$ is a cyclic solution to $(K/k,m,S)$ then $L_p$ clearly is a cyclic solution to $(K_p/k,m_p,S_p)$.
  Conversely, if for each $p|m$ there  is a cyclic solution $L^{(p)}$ to $(K_p/k,m_p,S_p)$
  then the field compositum of $K$ and all $L^{(p)}$ is a cyclic solution to $(K/k,m,S)$.
\end{proof}

\begin{lemma}
  Suppose that $\mu_m\subset k$.
  Then $\mu_m\subset N_{K/k}(K)$ if and only if $\mu_{m_p}\subset N_{K_p/k}(K_p)$ for each prime divisor $p$ of $m$.
  \label{lem:mp}
\end{lemma}
\begin{proof}
 The ``only if'' part is trivial. 
 The ``if'' part follows from two easy facts using induction on the number of prime divisors in $m$.
 First, if $\mu_{m_p}\subset N_{K_p/k}(K_p)$ then $\mu_m\subset N_{K_p/k}(K_p)$.
 Second, if $K_1/k$ and $K_2/k$ are two extensions of relatively prime degrees
 then $x\in N_{K_1/k}(K_1)\cap N_{K_2/k}(K_2)$ implies $x\in N_{K'/k}(K')$ where $K'$ is the field compositum $K_1K_2$.
\end{proof}

Next, we show in the following proposition that Theorem \ref{thm:cyclic}
further reduces to the case $S=\emptyset$
which is a global cyclic embedding problem without any local condition.

\begin{prop}
  \label{prop:cyclic}
  $(K/k,m,S)$ has a cyclic solution if and only if $(K/k,m,\emptyset)$ has.
\end{prop}
\begin{rem*}
  Proposition \ref{prop:cyclic} does not hold if one replaces
  ``cyclic solution'' with ``abelian solution'' or ``solution''.
  As an example, $(\Q(\zeta_7+\zeta_7^{-1})/\Q,3,\set{2})$ has no solution
  but $(\Q(\zeta_7+\zeta_7^{-1})/\Q,3,\emptyset)$ clearly has abelian solutions.
\end{rem*}
Before starting the proof we briefly recall some well-known facts about
the correspondence between cyclic extensions and characters of absolute Galois groups.
Let $k$ in this paragraph be any field,
$\wt k$ its separable closure
and $G=\Gal(\wt k/k)$ its absolute Galois group, a pro-finite group.
Denote by $X(G)$ the abelian group $\Hom_c(G,\C^*)$ of continuous characters of $G$
(where $\C^*$ is equipped with the discrete topology).
Any $\varphi\in X(G)$ defines an extension $k(\varphi)/k$ of degree $\ord\varphi$
by taking $k(\varphi)$ the fixed field of $\ker\varphi$.
For, we have $\Gal(k(\varphi)/k)=G/\ker\varphi\cong\im\varphi$
and the image of any $\varphi\in X(G)$ is the finite group $\mu_n(\C^*)$ where $n=\ord\varphi$.
Conversely, any cyclic finite extension of $k$ is obtained in this way.
Obviously, for any $l\in\N$, we have $k(\varphi^l)\subseteq k(\varphi)$ and 
\begin{equation} [k(\varphi):k(\varphi^l)]=(l,\ord\varphi). 
  \label{eq:subf}
\end{equation}
If $\varphi,\varphi'\in X(G)$ both have $p$-power order, $p$ prime,
and $\ord\varphi\neq\ord\varphi'$ then
\begin{equation} \ord(\varphi\varphi')=\max\set{\ord\varphi,\ord\varphi'}.
  \label{eq:ord-max}
\end{equation}
For any valuation $v$ on $k$ denote by $G_v$ the absolute Galois group of $k_v$
and by $\varphi_v$ the restriction of $\varphi\in X(G)$ to the subgroup $G_v$.
If $w$ is a valuation on $k(\varphi)$ and $w|v$ then
\begin{equation}
  k(\varphi)_w\cong k_v(\varphi_v) 
  \label{eq:v}
\end{equation}

\begin{proof}[Proof of Proposition \ref{prop:cyclic}]
  (I thank Jack Sonn for pointing out this argument to me.)
  In view of Lemma \ref{lem:p} we assume that $m$ and $[K:k]$ are both powers of the same prime $p$.
  Let $m=p^s$. Trivially, we can assume $s>0$.
  We further assume without loss of generality that $S_0$ is nonempty (for a technical reason).

  Suppose $(K/k,p^s,\emptyset)$ has a cyclic solution $L$.
  Let $L=k(\psi)$ with $\psi\in X(G)$.
  Then $\ord\psi$ is a $p$-power $\geq p^s$.
  We conclude from \eqref{eq:subf} that $K=k(\varphi)$ for $\varphi:=\psi^{p^s}$.
  Let $w|v$ be valuations on $L$ and $K$ respectively.
  By \eqref{eq:v} we have $K_v\cong k_v(\varphi_v)$ and $L_w\cong k_v(\psi_v)$.
  Since $(\psi_v)^{p^s}=\varphi_v$, 
  \eqref{eq:subf} implies $[L:K]_w=(p^s,\ord\psi_v)$.
  Hence,
  \begin{equation*}
      \text{$L$ has full local degree in $S \iff \ord\psi_v\geq m_v$ for all $v\in S$}
  \end{equation*}
  where $m_v$ is defined $p^s$ if $v\in S_0$ and $(2,p)$ if $v\in S_\ift$.
  Let $S'$ be the subset of $S$ where this condition fails, i.e.
  $S':=\sett{v\in S}{\ord\psi_v<m_v}.$ 
  By the weak Grunwald-Wang Theorem
  (cf.\ Hasse \cite[Ganz schwacher Existenzsatz]{hasse:grunwald} or Artin-Tate \cite[Ch.\ 10, Thm.\ ~5]{artin-tate})
  there is a character $\chi\in X(G)$ of order $p^s$
  such that $\ord\chi_v=m_v$ for each $v\in S'$
  and $\ord\chi_v=1$ for each $v\in S\bs S'$.
  For each $v\in S$, since
  $\ord\psi_v\neq\ord\chi_v$ (unless $m_v=1$) 
  and $\max\set{\ord\psi_v,\ord\chi_v}\geq m_v$
  we have $\ord(\psi\chi)_v\geq m_v$ by \eqref{eq:ord-max}. 
  In particular, $\ord(\psi\chi)\geq p^s$ because $S_0$ was assumed to be nonempty.
  Clearly, $(\psi\chi)^{p^s}=\psi^{p^s}=\varphi$, thus $k(\psi\chi)$ is a solution to $(K/k,p^s,S)$. 
\end{proof}

The cyclic embedding problem is classically answered over arbitrary fields by the following two theorems.
Let $K/k$ be any finite cyclic extension and let $m\in\N$.

\begin{theorem}[Albert's Theorem]
  Suppose $\charak k\ndiv m$ and $\mu_m\subset k$.
  Then $K/k$ embeds into a cyclic extension $L/k$ with $[L:K]=m$ 
  if and only if $\mu_m\subset N_{K/k}(K)$ and $k$ possesses a cyclic extension of degree $m_0$,
  where $m_0$ is the maximal divisor of $m$ relatively prime to $[K:k]$.
  \label{thm:albert}
\end{theorem}

\begin{rem*}
Theorem \ref{thm:albert} is due to Albert \cite{albert:cyclic-fields} with $m$ prime
and the proof is reproduced in Albert's book \cite[Ch.\ 9, Thm.\ 11]{albert:modern-algebra}.
For convenience to the reader a modification of this proof to general $m$ is placed in the appendix (\S \ref{sec:albert}).
A cohomological proof can be found in Arason-Fein-Schacher-Sonn \cite[Prop.\ 1]{afss:height}.
\end{rem*}

\begin{theorem}[Artin-Schreier-Albert-Witt]
  Suppose $\charak k=p$ and $m=p^r$ with $r>0$.
  Then $K/k$ embeds into a cyclic extension $L/k$ with $[L:K]=m$ 
  if and only if $k$ possesses a cyclic extension of degree $p$.
  \label{thm:asaw}
\end{theorem}

\begin{rem*}
Theorem \ref{thm:asaw} is due to Albert \cite{albert:cyclic-p} 
with an inductive argument based on Artin-Schreier theory.
The same kind of argument is found in Witt \cite[\S\ V]{witt:konstr-gal-koerpen}
or in the books of Albert \cite[p.\ 194f]{albert:modern-algebra} 
and Jacobson \cite[Thm.\ 4.2.3, p.\ 159]{jacobson:fin-dim-div-alg}.
Alternatively, the result is obtained as an immediate consequence of the theory of Witt vectors
without the necessity of induction
(see e.g.\ Jacobson \cite[\S\ 8.11]{baii}).
\end{rem*}

\begin{proof}[Proof of Theorem \ref{thm:cyclic}]
As a corollary of Theorem \ref{thm:asaw} we can drop the assumption $\charak k\ndiv m$ in Theorem~\ref{thm:albert}.
To see this, suppose $\charak k=p$ and $m=p^rm'$ with $p\ndiv m'$.
Since $\mu_m(k)=\mu_{m'}(k)$, the ``only if''-part is trivial.
Conversely, let $\mu_{m'}\subset N_{K/k}(K)$ and suppose $k$ possesses a cyclic extension of degree $m_0$, where $m_0$ is the maximal divisor of $m$ relatively prime to $[K:k]$.
By Theorem \ref{thm:albert}, $K/k$ embeds into a cyclic extension $L'/k$ with $[L':K]=m'$.
Applying Theorem \ref{thm:asaw} to $L'/k$ proves the claim, we only need to verify that $k$ possesses a cyclic extension of degree $p$.
If $p$ divides $[K:k]$ then such extension is contained in $K$.
And if $p\ndiv[K:k]$ then $p|m_0$, hence a cyclic extension of degree $p$
exists by hypothesis.

Together with Proposition \ref{prop:cyclic} this completes the proof of Theorem \ref{thm:cyclic}.
Indeed, a global field $k$ has cyclic extensions of any degree.
Hence, Theorem \ref{thm:albert} (without the assumption on $\charak k$)
states that $(K/k,m,\emptyset)$ has a cyclic solution if and only if $\mu_m\subset N_{K/k}(K)$.
\end{proof}

\begin{rem*}
  Proposition~\ref{prop:cyclic} --- on which Theorem \ref{thm:cyclic} relies ---
  makes use of the weak Grunwald-Wang Theorem.
  However, since Theorem \ref{thm:cyclic} assumes $\mu_m\subset k$,
  the weak Grunwald-Wang Theorem needs to be invoked within the proof of Proposition \ref{prop:cyclic}
  only for degrees $p^s$ with $\mu_{p^s}\subset k$.
  It is well-known that this case of the Grunwald-Wang Theorem (even the strong version) 
  can be proved quite elementally with Kummer theory
  (cf.\ Hasse \cite[\S 4]{hasse:grunwald}).
  Therefore, the claim of the introduction 
  that the Main Theorem requires only Kummer theory plus the approximation theorem is not affected
  by this use of Grunwald-Wang.
\end{rem*}

\begin{rem*}
  The case of Theorem \ref{thm:cyclic} where $m$ is a power of $\charak k$ was originally proved in \cite[Thm.\ 6.7]{hanke:thesis}.
\end{rem*}

\subsection{General and abelian solutions}
Let $K/k$ be a cyclic extension of global fields and let $m\in\N$.

\begin{theorem}
For any $S$,
\begin{enumerate}
  \item 
 if $\mu_m\subset K$ then $(K/k,m,S)$ has a solution,
 \item
 if $\mu_m\subset k$ then $(K/k,m,S)$ has an abelian solution.
\end{enumerate}
  \label{thm:main2}
\end{theorem}
\begin{proof}
By taking field composita the proof reduces to the case $m=p^s$, $p$ prime. 
We assume $p\neq\charak k$, 
for otherwise there is even a cyclic solution by Theorem \ref{thm:cyclic}.

Suppose $\mu_{p^s}\subset K$.
By Theorem \ref{thm:exist-a} and Corollary \ref{cor:ass}, we can choose for each $v\in S_0$ 
an element $a_v\in k_v$ such that $x^{p^s}-a_v$ is irreducible over $K_v$.
For each $v\in S_\ift$ choose $a_v=-1$.
Recall that in a complete valued field $k$ with $\charak k\nmid n$ 
any element sufficiently close to $1$ is an $n$-th power\footnote{
See Efrat \cite[Prop.\ 18.2.1]{efrat:val} for a proof.
Note that this is trivial for archimedian valuations.}.
Using the approximation theorem, 
let $a\in k$ such that $\frac{a}{a_v}=1+\frac{a-a_v}{a_v}$ 
is a $p^s$-th power in $k_v$ for each $v\in S$.
By Kummer theory we have $K_v(a^{1/{p^s}})=K_v(a_v^{1/{p^s}})$.
Assuming without loss of generality that $S_0$ is nonempty
it follows that $L:=K(a^{1/{p^s}})$ has degree $p^s$ over $K$ and full local degree in $S$.
Since $K/k$ is Galois and $a\in k$ we have $L/k$ Galois,
hence $L$ is a solution to $(K/k,p^s,S)$. 
Moreover, if $\mu_{p^s}\subset k$ then $L/k$ is abelian because $L$ is the compositum of the two cyclic extensions $K/k$ and $k(a^{1/p^s})/k$.
Thus (i) and (ii) are proved.
\end{proof}

\section{Application to central simple algebras}

The Main Theorem translates as follows into a statement about absolute Galois splitting fields of central simple algebras.

\begin{theorem}
  Let $K/k$ be a cyclic extension of global fields and let $A$ be a central simple $K$-algebra.
  If $A$ is a tensor product over $K$ of a symbol algebra with a $p$-algebra
  (the factors are allowed to be trivial)
  then $A$ contains a $k$-Galois maximal subfield.
  \label{thm:main-cs}
\end{theorem}
\begin{proof}
  Let $m=\deg A$ and define $S$ to be the set of normalized valuations on $K$ for which the Hasse invariant is nonzero.
  Then $S$ is a finite set (cf.\ Pierce \cite[\S 18.5]{pierce:ass-alg})
  of pairwise inequivalent valuations
  that contains no complex valuations (Pierce \cite[\S 18.5]{pierce:ass-alg}).
  If $A$ is as assumed then obviously $\mu_m\subset K$,
  hence $(K/k,m,S)$ has a solution by Theorem \ref{thm:main2} (i).
  Any solution $L$ to $(K/k,m,S)$ is a splitting field of $A$ by the local-global principle
  (cf.\ \cite[\S18.4]{pierce:ass-alg}).
  Moreover, $L$ is isomorphic to a maximal subfield of $A$ because $[L:K]=\deg A$.
  This proves the Theorem.
\end{proof}

\begin{rem*}
  In fact, Theorem \ref{thm:main-cs} is equivalent to Theorem \ref{thm:main2} (i).
  For, if $K$ is a global field and $\mu_m\subset K$ 
  then any central simple $K$-algebra of degree $m$ is a 
  tensor product of a symbol algebra with a $p$-algebra.
  Furthermore, for any $S$ there exists a central simple $K$-algebra of degree $m$ 
  that has a Hasse invariant of order $m$ at all $v\in S_0$ 
  and of order $(2,m)$ at all $v\in S_\ift$. 
  Theorem \ref{thm:main-cs} thus implies Theorem \ref{thm:main2} (i) by the local-global principle.
\end{rem*}

We first apply the result to twisted Laurent series rings and twisted function fields.
A brief recollection of these constructions follows;
for a thorough reference see Cohn \cite[\S 5.2]{cohn:intro-ring-th}
(other sources are Pierce \cite[\S 19.7]{pierce:ass-alg} and Jacobson \cite[\S 1.10]{jacobson:fin-dim-div-alg} but they treat only Laurent series). 
Let $A$ be a central simple $K$-algebra and let $\sigma$ be an automorphism of $A$
whose restriction $\sigma|_K$ is of finite order.
The twisted Laurent series ring $A((\bf x;\sigma))$ is the ring of Laurent series over $A$
where the indeterminate ${\bf x}$ acts on $A$ via $\sigma$ : ${\bf x}d{\bf x}^{-1}=\sigma(d)$.
A subring of $A((\bf x;\sigma))$ is formed by the twisted polynomials,
and the ring of central quotients of this subring is called the twisted function field $A(\bf x;\sigma)$.
The centres of  $A((\bf x;\sigma))$ and $A(\bf x;\sigma)$ can be identified with $k((\bf t))$ and $k(\bf t)$, respectively,
where $k$ is the fixed field of $\sigma|_K$ and $\bf t$ is a commuting indeterminate.
Over these fields $A((\bf x;\sigma))$ and $A(\bf x;\sigma)$ are central simple algebras of degree $[K:k]\cdot\deg A$. 

\begin{cor}
  \label{cor:symbol2}
  If $K$ is a global field and $A$ is a tensor product over $K$ of a symbol algebra with a $p$-algebra then 
  $A( (\bf x;\sigma) )$ and $A(\bf x;\sigma)$ are crossed products.
\end{cor}
\begin{proof}
  Let $L$ be a $k$-Galois maximal subfield of $A$ by Theorem \ref{thm:main-cs}.
  Since $A$ is identified with $A{\bf x}^0$, $L$ is also a subfield of $A((\bf x;\sigma))$ and $A(\bf x;\sigma)$.
  Obviously, $L(({\bf t}))/k(({\bf t}))$ and $L({\bf t})/k({\bf t})$ are Galois extensions of degree $[L:k]$.
  Since $A((\bf x;\sigma))$ and $A(\bf x;\sigma)$ have degree $[K:k]\cdot\deg A=[L:k]$ these two algebras are crossed products.
\end{proof}

\begin{rem*}
  Corollary \ref{cor:symbol2} is no longer true for 
  iterated twisted Laurent series rings or function fields.
  Indeed, \cite{hanke:expl-ex} exhibits a noncrossed product example of the form $D((\bf x;\sigma))((\bf y;\tau))$ resp.\ $D(\bf x;\sigma)(\bf y;\tau)$ 
  with $D$ a quaternion division algebra.
\end{rem*}

\begin{rem*}
Corollary \ref{cor:symbol2} should be seen in contrast to the example 
from \cite{hanke:laurent-noncr}
of a noncrossed product twisted Laurent series ring $D((\bf x;\sigma))$
with $D$ a cyclic cubic division algebra.
A considerable difficulty in making the example from \cite{hanke:laurent-noncr} explicit
was the computation of the automorphism $\sigma$ of $D$,
achieved only with the help of a computer algebra system.
Since this task is easier for symbol algebras (cf.\ Hanke \cite[Prop.\ 3.2]{hanke:expl-ex})
it is natural to ask whether there are noncrossed products of the form $D( (\bf x;\sigma) )$ with $D$ a symbol algebra.
Corollary \ref{cor:symbol2} states that this is not the case (over global fields).
\end{rem*}

For division algebras Corollary \ref{cor:symbol2} can be viewed more generally in the context of valued division algebras
that are inertially split.
See Schilling \cite{schilling:th-of-val} for an introduction to valuations on division algebras
and Jacob-Wadsworth \cite[\S 5]{jacob-wadsworth:div-alg-hensel-field} for the theory of inertially split division algebras.
Let $k$ be a number field (we need $k$ to be perfect)
and let $(F,v)$ be a field with discrete Henselian rank $1$ valuation and residue field $k$.
Since $v$ is Henselian it extends to a valuation on any $F$-central division algebra $D$
(\cite[p.\ 53, Thm.\ 9]{schilling:th-of-val}).
Let $\ovl D$ denote the residue division algebra with respect to this extension. 

\begin{cor}
  \label{cor:hensel}
  Let $D$ be an $F$-central division algebra and let $m=\deg\ovl D$.
  If $\mu_m\subset Z(\ovl D)$ then $D$ is a crossed product.
  Moreover, if $\mu_m\subset k$ then $D$ is an abelian crossed product,
  and if $\mu_m\subset N_{Z(\ovl D)/k}(Z(\ovl D))$ then $D$ is cyclic.
\end{cor}
\begin{proof}
  Since $v$ is Henselian with respect to a discrete rank $1$ valuation and has perfect residue field
  every $F$-central division algebra is inertially split 
  (see Serre \cite[Ch.\ X, \S 7, Example~b)]{serre:local-fields}).
  The extension $Z(\ovl D)/k$ is cyclic (\cite[Lem.\ 5.1]{jacob-wadsworth:div-alg-hensel-field})
  and we have $[L:k]=m \cdot [Z(\ovl D):k]=\deg D$
  (\cite[Thm.\ 5.15]{jacob-wadsworth:div-alg-hensel-field}).
  
  Suppose $\mu_m\subset Z(\ovl D)$
  and let $L$ be a $k$-Galois maximal subfield of $\ovl D$ by Theorem~\ref{thm:main-cs}.
  Let $M$ be an inertial lift of $L$ over $F$ in $D$ (\cite[Thm.\ 2.9]{jacob-wadsworth:div-alg-hensel-field}).
  We have $[M:F]=[L:k]=\deg D$ and $M/F$ is Galois, hence $D$ is a crossed product.

  If $\mu_m\subset k$ then $L/k$ can be found abelian by Theorem \ref{thm:main2} (ii),
  and if $\mu_m\subset N_{Z(\ovl D)/k}(Z(\ovl D))$ then $L/k$ can be found cyclic by Theorem \ref{thm:cyclic}.
  The additional statements are immediate because $M/F$ and $L/k$ have isomorphic Galois groups.
\end{proof}

\section{Appendix}
\label{sec:albert}

Let $K/k$ be an arbitrary cyclic field extension and let $m\in\N$.

\begin{varthm}[Albert's Theorem]
  Suppose $\charak k\ndiv m$ and $\mu_m\subset k$.
  Then $K/k$ embeds into a cyclic extension $L/k$ with $[L:K]=m$ 
  if and only if $\mu_m\subset N_{K/k}(K)$ and $k$ possesses a cyclic extension of degree $m_0$,
  where $m_0$ is the maximal divisor of $m$ relatively prime to $[K:k]$.
  \label{thm:albert2}
\end{varthm}
\begin{proof}
  For $m$ prime the theorem is proved in Albert \cite[Ch.\ 9, Thm.\ 11]{albert:modern-algebra}.
  A cohomological proof can be found in Arason-Fein-Schacher-Sonn \cite[Prop.\ 1]{afss:height}.
  The following proof is along the lines of Albert's original proof
  and explicitly constructs the field $L$.

  Let $[K:k]=n$.
  If $L$ is as in the theorem then $L/K$ is a Kummer extension of degree $m$.
  Let $L=K(\alpha), \alpha^m\in K,$ and let $\sigma$ be a generator of $\Gal(L/k)$.
  Clearly, $\sigma^n(\alpha)=\zeta_m\alpha$ with $\zeta_m$ a primitive $m$-th root of $1$.
  It follows easily from the hypothesis $\mu_m\subset k$ that $a:=\sigma(\alpha)/\alpha$ is fixed by $\sigma^n$,
  hence $a$ lies in $K$. 
  Moreover, $N_{K/k}(a)=\zeta_m$.
  Since $L/k$ trivially contains a subextension of degree $m_0$ the ``only if'' part is proved.
  
  Conversely, suppose $\mu_m\subset N_{K/k}(K)$.
  By taking field composita we are reduced to the case $m_0=1$,
  i.e.\ we assume that every prime divisor of $m$ divides $[K:k]$.
  Let $a\in K$ such that $N_{K/k}(a)$ is a primitive $m$-th root of unity.
  By Hilbert's Theorem~90 there is $b\in K^*$ with $\sigma(b)/b=a^m$
  where $\sigma$ is a generator of $\Gal(K/k)$.
  We claim that the order of $b\in K^*$ modulo ${K^*}^m$ is $m$.
  Assume this is not the case, i.e.\ $b^{m/p}\in {K^*}^m$ for some $p|m$.
  The hypothesis $\mu_m\subset k$ implies that for any $d|m$ the homomorphism
  $(\cdot)^{d}:K^*/{K^*}^{m/d}\lra K^*/{K^*}^m$ is injective.
  Thus, $b\in {K^*}^p$.
  If $b=c^p$ with $c\in K^*$
  then $a^m=(\sigma(c)/{c})^p$, hence $a^{m/p}=\zeta\sigma(c)/{c}$
  for some $p$-th root of unity $\zeta\in k$.
  We derive $N_{K/k}(a)^{m/p}=\zeta^{n}=1$ because $p|n$.
  Since this contradicts the choice of $a$, the claim is proved.
  Thus, $L:=K(b^{1/m})$ has degree $[L:K]=m$. 
  To see that $L$ is cyclic over $k$ it is easily verified that mapping
  $b^{1/m}\lms ab^{1/m}$ defines an extension of $\sigma$ to $L$ of order $mn$.
\end{proof}

\bibliographystyle{plain}
\bibliography{master}

\end{document}